\documentclass{article}
\usepackage{amssymb}
\usepackage{amsmath}
\def\ul{\underline}
\def\ov{\overline}

\begin{document}
\centerline{\Large{\scshape Integral closures of ideals}}
\centerline{\scshape Douglas A. Leonard}
\centerline{Department of Mathematics and Statistics}
\centerline{Auburn University}
\centerline{Auburn, AL 36849-5307, USA}
\medskip
\begin{abstract}
The Qth-power algorithm 
for computing structured global presentations 
of integral closures of affine domains over finite fields
is modified to compute structured presentations 
of integral closures of ideals in affine domains over finite fields
relative to a local monomial ordering.
A non-homogeneous version of the standard (homogeneous) Rees algebra
is introduced as well. 
\end{abstract}

\section{Introduction}
There are two recent books by Swanson and Huneke \cite{Swan} and Vasconcelos \cite{Vas06} 
dealing with integral closures, with emphasis on integral closures of ideals.
There are also two implementations, 
{\tt normalI} in {\sc Singular}'s {\tt reesclos.lib} \cite{GPS}, \cite{Hir}  and 
{\tt integralClosure} in {\sc Macaulay2}'s {\tt IntegralClosure} package \cite{M2}, \cite{ic}.

The emphasis here will be on producing a version of the Qth-power algorithm 
\cite{qth}, \cite{Pell}, \cite{mod}, \cite{Qth}, \cite{char0},
that will produce the integral closure of an ideal in an affine domain over a finite field 
with structure induced from the original ring viewed over a fixed Noether normalization. 
(As with the Qth-power algorithm for ring integral extensions,
the restriction to finite fields can be relaxed \cite{char0}, 
but that will not be part of this paper.)
There is code written by the author for {\sc Macaulay2} implementing this approach.

Standard methods use the {\em Rees algebra} of the ideal 
to change an integral closure of ideals problem
into an integral closure of rings problem,
extracting the integral closure of the ideal from the integral closure of its Rees algebra.
If $C(A,Q(A))$ denotes the integral closure of the ring $A$ in its field of fractions,
then this approach actually produces the integral closures $C(I^k,C(A,Q(A)))$  
of the ideals $I^k$ of the integral closure $C(A,Q(A))$ of the ring $A$,
extracting the integral closures $C(I^k,A)$ as $C(I^k,C(A,Q(A)))\cap A$.

Since $C(I^k,C(A,Q(A)))$ would seem to be a more interesting object than $C(I^k,A)$,
and is the object naturally produced, it is what shall be focused on here. 
A Gr\"obner basis for a larger ring 
(called rees(I) here for want of a better name) 
than the Rees algebra will be computed.
And since the Qth-power algorithm deals with $P$-module presentations 
for $P$ a Noether normalization of $A$,
the input and output will have this flavor.

While computing a presentation of an integral closure of a ring
makes better sense using global monomial orderings, 
local monomial orderings are more suitable 
for computing the integral closure of an ideal,
as the structure of Rees algebras might suggest.

The approach used to compute integral closures of ideals here 
will be to start with $\ov{A}$ integrally closed (or compute $\ov{A}:=C(A,Q(A))$)
using a local presentation,
then compute a nested sequence of ideals 
$I^{k-1}=:I^k_0\subseteq\cdots\subseteq I^k_L=I^k_{L+1}=C(I^k,\ov{A})$ for each $k$
as long as $C(I^{k},\ov{A})\neq C(I^{k-1},\ov{A})C(I,\ov{A})$.
Each of these ideals will be produced by computing sequences of nested $P$-modules
 \[ I^{k-1}=:M_0(I^k_j)\supseteq\cdots\supseteq M_L(I^k_j)=M_{L+1}(I^k_j)=:I^k_{j+1},\]
with
\[ M_{i+1}(I^k_j):=\{g\in M_{i}(I^k_j)\ :\ g^Q\in M_{i}(I^k_j)^{Q-1}I^k_j\}.\]
Both sequences are means-to-an-end, meaning they seem to have little independent import.

\section{Integrality}
Start with a polynomial ring $P:=\mathbf{F}_q[x_n,\ldots,x_1]$ 
in $n$ (independent) variables over a finite field $\mathbf{F}_q$ of prime order $q$. 
Adjoin new (dependent) variables to get an integral extension $A:=R/J$ 
with $R:=P[y_m,\ldots,y_1]$ 
(written as the {\em flattened polynomial ring} $\mathbf{F}_q[\ul{y};\ul{x}]$)  
and $J$ an ideal of $R$ 
such that each $y\in A$ is {\em integral over P}, 
meaning it satisfies some monic polynomial
$f_y(T):=T^k+a_1T^{k-1}+\cdots+a_kT^0\in P[T]$ of some degree $k$.

The {\em integral closure} of $A$ in its {\em field of fractions} $Q(A)$, denoted $C(A,Q(A))$
is the ring of all elements of $Q(A)$ integral over $A$. 
Assume from now on that $A=C(A,Q(A))$ is integrally closed, but write it as $\ov{A}$
to remind ourselves of this assumption.
[An ideal $I$ of $A$ can be thought of as an ideal of $\ov{A}$
(though technically it should probably be written as $I\ov{A}$ 
to reflect that it is closed under multiplication by elements of $\ov{A}$ 
and not just by elements of $A$.]

Let $I:=\langle g_0,\ldots,g_s\rangle$ be an ideal if $\ov{A}$.
$z$ is {\em integral over I} iff it satisfies some monic {\em I-polynomial}
$f_z(T):=T^k+a_1T^{k-1}+\cdots+a_kT^0\in \ov{A}[T]$ of some degree $k$, 
with the more restrictive condition that $a_j\in I^j$ for $1\leq j\leq k$.

\section{Non-homogeneous Rees algebras}

The standard (homogeneous) {\em Rees algebra}, Rees(I), 
often written as $\ov{A}[It]$ or $\ov{A}[g_0t,\ldots,g_st]$ or $\oplus_{k\geq 0}I^kt^k$,
is the image of the map $\phi\ :\ \ov{A}[G_0,\ldots,G_s]\rightarrow \ov{A}[t]$ 
defined by $\phi(G_j):=g_jt$ for $0\leq j\leq s$. 
[The notation $\ov{A}[It]$ is just shorthand for $\ov{A}[g_0t,\ldots,g_st]$.
$\ov{A}[g_0t,\ldots,g_st]$ itself is problematic in practice 
in that $g_jt$ is not necessarily an acceptable name for a variable. 
And $\oplus_{k\geq 0}I^kt^k$ is also problematic 
in that elements such as $(x_{3,1}^2)^3(x_{2,0}^3)^2=(x_{3,1}x_{2,0})^6$ 
in an example below will correspond to, among other things, 
$(x_{3,1}^2t)^3(x_{2,0}^3t)^2\in I^5t^5$ and $(x_{3,1}x_{2,0}t)^6\in I^6t^6$.
Even using the {\em associated graded ring} $R/I\oplus I/I^2\oplus\cdots $ 
doesn't necessarily change this.] 

The corresponding {\em presentation} of Rees(I) 
is $\ov{A}[G_0,\ldots,G_s]/ker(\phi)$.
[This doesn't resolve the question above.
Having the homogeneous relation $G_{6,2}^3G_{6,0}^2-G_{5,1}^5x_{3,1}x_{2,0}\in ker(\phi)$
is not the same as having  $G_{6,2}^3G_{6,0}^2\mapsto G_{5,1}^6$.]
But the important observation here is that the map $\phi$ is ``external'' to this presentation 
in that for any $0\neq \alpha\in \ov{A}$ and $\phi_{\alpha}(G_j):=\alpha g_jt$, 
the presentation will be the same,
since a minimal, reduced Gr\"obner basis for $ker(\phi)$ 
necessarily consists of elements that are {\em homogeneous} in the variables $G_0,\ldots,G_s$.

Consider instead a larger (non-homogeneous) quotient ring 
\[rees(I):=\ov{A}[G_0,\ldots,G_s][t^{-1}]/\langle g_j-G_jt^{-1}\rangle\]
in which the map $\phi$ is ``internal'' 
as it is defined by the generator relations of the quotient ring.
[The use of $t^{-1}$ instead of $t$ is supposed to suggest a local ordering 
in which $1\succ t^{-1}\succ\cdots$;
but in either case $t$ or $t^{-1}$ could be supressed 
if elements are implicitly graded appropriately.]
And we shall find a sequence of overrings of this, 
terminating in a {\em non-homogeneous} Rees algebra, $rees(C(I,\ov{A}))$. 

{\bf Definition 1} Given an ideal
$I:=\langle g_i\ :\ 0\leq i\leq s\rangle\in \ov{A}$,
recursively define $rees_k(I,\ov{A})$,
starting with $rees(I^0,\ov{A})$ a local presentation of $\ov{A}$, as follows.
If $\{ g_{i,k}\ :\ 0\leq i\leq s(k)\}$ is a generating set
for those elements of $C(I^k,\ov{A})$ not already in $C(I^{k-1},\ov{A})C(I,\ov{A})$,
then 
\[rees_k(C(I,\ov{A})):=rees_{k-1}(C(I,\ov{A}))[ G_{i,k}\ :\ 0\leq i\leq s(k)]/
\langle g_{i,k}t^{1-k}-G_{i,k}t^{-k}\ :\ 0\leq i\leq s(k)\rangle.\]
Then $rees(C(I,\ov{A}))$ is the (finite) union of the $rees_k(C(I,\ov{A}))$.

[Note that while Rees algebras have the appealing theoretical property that 
$Rees(C(I,\ov{A}))=\sum_kC(I^k,\ov{A})t^{-k}$ 
(a generating function for the sequence $C(I^k,\ov{A})$),
$rees(C(I,\ov{A}))$ has the following advantage, 
as noted in the examples above and reiterated in the final example below.]

{\bf Theorem 2} If $f\in\ov{A}$, then $f\in C(I^k,\ov{A})$ 
iff $t^{-k}$ divides $NormalForm(f,J(I,\ov{A}))$ 
for $J(I,\ov{A})$ the ideal of relations defining $rees(C(I,\ov{A}))$.

{\bf Proof :} Any $f\in\ov{A}$ is clearly reduced to  $NormalForm(f,J(I,\ov{A}))$ 
explicitly an element of  $C(I^k,\ov{A})$ for largest possible $k$. 

It would also seem that the closures $C(I^k,\ov{A})$ are more closely related to
local rather than global presentations of $\ov{A}$, 
given the power series flavor of Rees algebras 
and the desire to write an element as living in the highest power $I^k$ possible.
A simple example such as $I:=\langle x^5+x^2,y^5+y^2\rangle\in \mathbf{F}[x,y]$
shows that the integral closure is related to that of
 $\langle x^2,y^2\rangle$, not that of $\langle x^5,y^5\rangle$.
So it is better to write $I=\langle x^2u^{-1},y^2v^{-1}\rangle=\langle x^2,y^2\rangle$
with $u^{-1}:=1+x^3$ and $v^{-1}:=1+y^3$ being {\em units}.
rees(I) is generated by $x^2-G_{2,2}t^{-1}$ and $y^2-G_{2,0}t^{-1}$, 
with $(xy)^2\mapsto G_{2,2}G_{2,0}t^{-2}$ showing that $xy\in C(I,\ov{A})$.
So $rees(C(I,\ov{A}))$ has additional relations
$xy-G_{2,1}t^{-1}$, and $(G_{2,1}^2-G_{2,2}G_{2,0})t^{-2}$.
[As a warning, it should be noted that local monomial orderings are different 
in flavor from global monomial orderings. 
A simple example such as $I:=\langle 1+x^2,1+y^2\rangle\in \mathbf{F}_2[x,y]$ 
will reduce locally to $I=\langle 1\rangle$ since both generators will be units. 
Homogenizing this to $I:=\langle h^2+x^2,h^2+y^2\rangle$ 
will produce the expected integral closure 
$C(I,\ov{A})=\langle h^2+x^2,h^2+hy+hx+yx,h^2+y^2\rangle$.]

\section{The Qth-power algorithm for ideals}

In this context, the Qth-power algorithm  
(using exponent $Q=q^e$ for sufficiently large $e$ in characteristic $q$),
naturally produces a nested sequence of $P$-modules 
 \[ C(I^{k-1},\ov{A})=:M_0(I^k_j)\supseteq\cdots\supseteq M_L(I^k_j)=M_{L+1}(I^k_j)=:I^k_{j+1},\]
with
\[ M_{i+1}(I^k_j):=\{g\in M_{i}(I^k_j)\ :\ g^Q\in M_{i}(I^k_j)^{Q-1}I^k_j\}\]
as mentioned above.

This requires only the computation of normal forms modulo $J$ and modulo $M_{i}^{Q-1}I$.
This is a $q$-linear algorithm in characteristic $q$.
And the normal forms require only a $P$-module generating set and Gr\"obner bases
based on the Gr\"obner basis for $rees(I,\ov{A})$.

That is, if $M_{i}=\langle b_0,\ldots,b_S\rangle$,
consider $\ov{A}[G_0,\ldots,G_s;B_0,\ldots,B_S]$
with ideal 
\[J_{q,i}:=\langle b_0-t^{-1}B_0,\ldots,b_S-t^{-1}B_S\rangle
         +\langle g_0-t^{-1}G_0,\ldots,g_s-t^{-1}G_s\rangle.\]
Then for $g\in M_{i}$, compute $NF(g^Q,J_{i})$;
and put $g\in M_{i+1}$ iff $NF(g^Q,J_{i})=0$.

{\bf Theorem 3} 
\[ M_{0}(I):=\langle 1\rangle\supset\cdots\supset M_{L}(I)=M_{L+1}(I)=:\phi(I)\]
is a nested sequence of $P$-modules.

{\bf Proof :} The nesting is obvious from the definition. 
That it is finite is because the set of leading monomials 
of the $P$-module generators are all subsets of 
$\{ y^i\ul{x}^{\ul{\alpha}}\ :\ 0\leq i<d, \alpha_j\leq\beta_j\}$
for $\ul{x}^{\ul{\beta}}$ the max of such occurring 
in the leading monomials of the generators of $I$.
There are only finitely many distinct such subsets, 
with each distinct $P$-module corresponding to a different subset.

Note that the small example 
\[I:=\langle a^5b^4,b^5c^4,c^5a^4\rangle\subset A:=\mathbf{F}[a,b,c]\]
has elements $a^4bc^4$, $a^4b^4c^1$, $a^1b^4c^4$, $a^4b^2c^3$,
$a^3b^4c^2$, and $a^2b^3c^4$ satifying $T^{21}-\gamma=0$ for some $\gamma\in I^{21}$;
elements $a^4b^3c^2$, $a^2b^4c^3$, and $a^3b^2c^4$ satisfying $T^7-\gamma=0$
for some $\gamma\in I^7$; and
$a^3b^3c^3$ satisfying $T^3-(a^5b^4)(b^5c^4)(c^5a^4)=0$.

But the Qth-power algorithm with a too small exponent, $Q=q=2$, produces
\[ \phi_2(I)=\langle a^4b^4c,\ a^4b^3c^2,\ a^4b^2c^3,\ a^4bc^4,\ a^3b^4c^2,\ a^3b^2c^4,\ a^2b^4c^3,\ a^2b^3c^4,\ ab^4c^4\rangle +I\]
from
\[(a^4b^3c^2)^2=(a^5b^4)(a^2b^3c^4),\ (a^2b^4c^3)^2=(b^5c^4)(a^4b^2c^3),\ (a^3b^2c^4)^2=(c^5a^4)(a^3b^4c^2),\]
\[(a^4b^2c^3)^2=(c^5a^4)(a^4b^4c^1),\ (a^3b^4c^2)^2=(a^5b^4)(a^1b^4c^4),\ (a^2b^3c^4)^2=(b^5c^4)(a^4b^1c^4),\]
\[(a^4b^4c^1)^2=(a^5b^4)(a^3b^4c^2),\ (a^1b^4c^4)^2=(b^5c^4)(a^4b^2c^3),\ (a^4b^1c^4)^2=(a^4b^5)(a^3b^4c^2);\]
while \[ (\phi_2)^2(I)=\phi_2(I)+\langle a^3b^3c^3\rangle\]
from 
\[ (a^3b^3c^3)^2=(a^3b^2c^4)(a^3b^4c^2).\]
For larger values of $Q$ only $\phi_Q(I)$ may be necessary.

{\bf Theorem 4}

Let $\phi_Q(I)$ denote the ideal produced from $I$ by the Qth-power algorithm 
with exponent $Q:=q^e$ in characteristic $q$. 
Then 
\[I\subseteq \phi_Q(I)\subseteq\cdots\subseteq 
      (\phi_Q)^k(I)=(\phi_Q)^{k+1}(I)=C(I,\ov{A}).\]

{\bf Proof :} The inclusions are obvious. 
And an ascending chain of ideals stabilzes by the ascending chain condition.
So the only thing to prove is that it stabizes at $C(I,\ov{A})$.

Let $f_1,\ldots,f_s$ be a finite generating set for $C(I,\ov{A})$.
Let $f_i^{D_i}+a_1f_i^{D_i-1}+\cdots+a_{D_i}=0$, for $1\leq i\leq s$.
Let $Q\geq\max\{ D_i\ :\ 1\leq i\leq s\}$.
Then $f_i^Q=-f_i^{Q-D_i}(a_1f_i^{D_i-1}+\cdots+a_{D_i})\in C(I,\ov{A})^{Q-1}I$,
meaning that $\langle f_1,\ldots,f_s\rangle\in M_j$ for all $j$ 
when computing $(\phi_Q)^{k+1}(I)$.
So $C(I,\ov{A})\subseteq (\phi_Q)^{k+1}(I)$, with the other inclusion obvious.

[Note that while one might be tempted to conjecture that $Q=q$ 
suffices for all characteristics $q$, there is a simple monomial counterexample
$I:=\langle a_1^{M-1}a_2,\ a_2^{M-1}a_3,\ a_{M}^{M-1}a_1\rangle$, which takes $Q\geq M$, irrespective
of the characterisitic $q$ to get $(a_1\cdots a_{M})^{M}=(a_1^{M-1}a_2)\cdots(a_{M}^{M-1}a_1)$.
So making the easy choice of $Q=q$ is not guaranteed to produce the whole integral closure.]

\section{Weighted presentations}
For integral closures of rings problems, a global weighted monomial ordering
allows for a nice strict affine $P$-algebra presentation in which $wt(f)=wt(NormalForm(f))$. 
For integral closures of ideals, since a local monomial ordering seems more appropriate, 
and since $\ov{A}$ is called for as a starting point for the computation, 
it makes more sense to ask for a local presentation. 

The definition of a global weight function:

{\bf Definition 5} A function $wt\ :\ R\backslash \{0\}\rightarrow \mathbf{N}$ 
is a ({\em global}) {\em weight function} iff the following hold.
\begin{enumerate}
   \item For $c\in \mathbf{F}$, $wt(c)=0$.
   \item $wt(ab)=wt(a)+wt(b)$.
   \item If $a+b\neq 0$ then $wt(a+b)\leq \min\{wt(a),wt(b)$, with equality if $wt(a)\neq wt(b)$.
   \item If $wt(a)=wt(b)$ then there exists $c\in \mathbf{F}$ such that 
         either $a=cb$ or $wt(a-cb)<wt(a)$.
\end{enumerate}
can be redone for a local monomial ordering as 

{\bf Definition 6a} A function $wt\ :\ R\backslash \{0\}\rightarrow -\mathbf{N}$ 
is a {\em local weight function} iff the following hold.
\begin{enumerate}
   \item For $c\in \mathbf{F}$, $wt(c)=0$.
   \item $wt(ab)=wt(a)+wt(b)$.
   \item $wt(a+b)\leq \min\{wt(a),wt(b)$, with equality if $wt(a)\neq wt(b)$.
   \item If $wt(a)=wt(b)$ then there exists $c\in \mathbf{F}$ such that 
         either $a=cb$ or $wt(a-cb)<wt(a)$.
\end{enumerate}

Or, if negative numbers are to be avoided this can be reworded as follows.

{\bf Definition 6b} A function $wt\ :\ R\backslash \{0\}\rightarrow \mathbf{N}$ 
is a {\em local weight function} iff the following hold.
\begin{enumerate}
   \item For $c\in \mathbf{F}$, $wt(c)=0$.
   \item $wt(ab)=wt(a)+wt(b)$.
   \item If $a+b\neq 0$ then $wt(a+b)\geq \min\{wt(a),wt(b)$, with equality if $wt(a)\neq wt(b)$.
   \item If $wt(a)=wt(b)$ then there exists $c\in \mathbf{F}$ such that 
         either $a=cb$ or $wt(a-cb)>wt(a)$.
\end{enumerate}

In any case, while the final condition of weights is nice, 
it is not necessary for the theory here, but only for providing a nicer, weighted
presentation of the quotient ring $\ov{A}$.

\section{Example, showing the recursive levels, $rees_k(C(I,\ov{A}))$}
Consider the following generic small example. 
Let $P:=\mathbf{F}_2[x_2,x_1]$.
The monic polynomial $f(T):=T^2+Tx_2^4x_1+x_2^2x_1^{11}\in P[T]$
defines an integral extension 
with $R:=P[y]$, $J:=\langle f(y)\rangle$, and $A:=R/J$.
Viewed as 
\[A=\mathbf{F}_2[y;x_2,x_1]/\langle f(y)\rangle\]
this has a weight function defined by the matrix
\[W(A):=\begin{pmatrix}13&2&2\\2&2&0\end{pmatrix}\]
naturally extending the grevlex monomial ordering on $P$.

The integral closure of this ring contains $z:=y/(x_2x_1)$, 
so $y$ is unnecessary in the presentation 
\[ \ov{A}:=\mathbf{F}_2[z;x_2,x_1]/\langle z^2+x_1^9+zx_2^3\rangle\]
of $C(A,Q(A))$, but only for the inclusion map from $A$ into $C(A,Q(A))$.
[That is, there may be several choices for $A$ 
having the same presentation for their integral closure.]

This has induced local weight function given by 
$W(\ov{A}):=\begin{pmatrix}9&3&2\\0&1&0\end{pmatrix}$,
so we shall use the names $y_{9,0}$, $x_{3,1}$, and $x_{2,0}$ 
instead of $z$, $x_2$, and $x_1$.

The relation defining $\ov{A}$ is:
\[ y_{9,0}^2+x_{2,0}^9+y_{9,0}x_{3,1}^3.\]
The non-homogeneous relations defining rees(I) are:
\[ x_{3,1}x_{2,0}-G_{5,1},\ x_{2,0}^3-G_{6,0},\ x_{3,1}^2-G_{6,2},\]
inducing the further relations:
\[ G_{5,1}x_{3,1}-G_{6,2}x_{2,0},\ G_{6,2}G_{6,0}-G_{5,1}^2x_{2,0},\ G_{6,0}x_{3,1}-G_{5,1}x_{2,0}^2,\]
\[y_{9,0}^2+G_{6,2}y_{9,0}x_{3,1}+G_{6,0}^3.\]
The additional non-homogeneous relation defining $rees_1(C(I,\ov{A}))$ is:
\[ y_{9,0}-G_{9,0},\]
inducing the further relation:
\[ G_{9,0}^2+G_{9,0}G_{6,2}x_{3,1}+G_{6,0}^3.\]
The further non-homogeneous relations defining $rees_2(C(I,\ov{A}))$ are:
\[ G_{9,0}x_{3,1}-G_{12,1},\ G_{9,0}x_{2,0}^2-G_{13,0},\]
inducing the further relation:
\[G_{9,0}G_{5,1}-G_{12,1}x_{2,0},\ 
  G_{9,0}G_{6,0}-G_{13,0}x_{2,0},\ 
  G_{9,0}G_{6,2}-G_{12,1}x_{3,1},\ 
  G_{9,0}^2+G_{6,0}^3+G_{12,1}G_{6,2},\]
\[G_{12,1}G_{6,0}-G_{13,0}G_{5,1},\ 
  G_{12,1}G_{9,0}+G_{6,0}^2G_{5,1}x_{2,0}^2+G_{12,1}G_{6,2}x_{3,1},\ 
  G_{12,1}^2+G_{6,0}^2G_{5,1}^2x_{2,0}+G_{12,1}G_{6,2}^2,\]
\[G_{13,0}x_{3,1}-G_{12,1}x_{2,0}^2,\ 
  G_{13,0}G_{6,2}-G_{12,1}G_{5,1}x_{2,0},\ 
  G_{13,0}G_{9,0}+G_{6,0}^3x_{2,0}^2+G_{12,1}G_{5,1}^2,\] 
\[G_{13,0}G_{12,1}+G_{6,0}^3G_{5,1}x_{2,0}+G_{12,1}G_{6,2}G_{5,1}x_{2,0},\ 
  G_{13,0}^2+G_{6,0}^4x_{2,0}+G_{12,1}G_{5,1}^2x_{2,0}^2.\] 
 
\section{Example with a local ordering with units, and the use of the rees algebra}

A slightly more complicated example than the one in the introduction is:
\[ A:=\mathbf{F}[x,y,z]/\langle x^6+x^3z-y^3z^2\rangle,\ I:=\langle xy^2,z\rangle.\]
For a global presentation of the integral closure $\ov{A}$, rewrite it as
\[ A:=\mathbf{F}[y_{5,3};x_{6,6},z_{6,0}]/\langle y_{5,3}^6-x_{6,6}^3x_{6,0}^2+y_{5,3}^3x_{6,0}\rangle, 
I:=\langle y_{5,3}x_{6,6}^2,x_{6,0}\rangle\]
with subscripts corresponding to the global weights.
The integral closure has elements $y_{9,9}:=x_{5,3}^3/x_{6,0}$,
$y_{8,6}:=(y_{5,3}^4+y_{5,3}x_{6,0})/(x_{6,6}x_{6,0})$, and
$y_{7,3}:=(y_{5,3}^5+y_{5,3}^2x_{6,0})/(x_{6,6}^2x_{6,0})$, with induced presentation
\[\ov{A}=\mathbf{F}[y_{9,9},y_{8,6},y_{7,3},y_{5,3};x_{6,6},x_{6,0}]/\ov{J}\]
with the ideal of induced relations $\ov{J}$ having minimal, reduced Gr\"obner basis consisting of
\[y_{9,9}^2-x_{6,6}^3+y_{9,9},\]
\[y_{9,9}y_{8,6}-y_{5,3}x_{6,6}^2,\]
\[y_{9,9}y_{7,3}-y_{5,3}^2x_{6,6},\]
\[y_{9,9}y_{5,3}-y_{8,6}x_{6,6}+y_{5,3},\]
\[y_{8,6}^2-y_{5,3}^2x_{6,6}-f_{7,3},\]
\[y_{8,6}y_{7,3}-y_{9,9}x_{6,0}-x_{6,0},\]
\[y_{8,6}y_{5,3}-y_{8,6}x_{6,6}+y_{5,3},\]
\[y_{7,3}^2-y_{8,6}x_{6,6}-y_{8,6}x_{6,0}+y_{5,3},\]
\[y_{7,3}y_{5,3}-x_{6,6}x_{6,0},\]
\[y_{5,3}^3-y_{9,9}x_{6,0}.\]
Rees(I) would have the $y_{5,3}x_{6,6}^2-G_{17,15}$, $x_{6,0}-G_{6,0}$, 
and any further relations induced by them.

There is a local presentation of $\ov{A}$, with relations
\[ y_{5,3}(1+y_{9,9})-y_{8,6}x_{6,6},\]
\[y_{9,9}(1+y_{9,9})-x_{6,6}^3,\]
\[x_{6,0}(1+y_{9,9})^2-y_{8,6}^3,\]
\[(y_{7,3}-y_{8,6}^2)(1+y_{9,9})^2-y_{8,6}^2x_{6,6}^3\]
with the global weights having no local meaning.
To clean this up, use the local variables and weights $x_{1,0}:=y_{8,6}$,
$x_{1,1}:=x_{6,6}$ and the unit $u:=1/(1+y_{9,9})$.
Then $y_{5,3}\mapsto x_{1,1}x_{1,0}u$, $x_{6,0}\mapsto x_{1,0}^3u^2$,
$y_{7,3}\mapsto x_{1,0}^2+x_{1,1}^3x_{1,0}^2u^2$,
$y_{9,9}\mapsto x_{1,1}^3u$, and the presentation has only
one remaining relation 
\[ 1-u-u^2x_{1,1}^3\]
defining the unit (as a power series in the local variable $x_{1,1}$).
The ideal $I\mapsto\langle x_{1,1}^3x_{1,0}u,x_{1,0}^3u^2\rangle$,
so rees(I) would have relations
\[ x_{1,0}^3-G_{3,0},\ x_{1,1}^3x_{1,0}-G_{4,3},\ G_{3,0}x_{1,1}^3-G_{4,3}x_{1,0}^2.\]

$x_{1,1}^2x_{1,0}^2$ completes the integral closure, 
so $rees(C(I,\ov{A})$ would have relations:
\[ x_{1,0}^3-G_{3,0},\ x_{1,1}^2x_{1,0}^2-G_{4,2},\ x_{1,1}^3x_{1,0}-G_{4,3},\]
\[ G_{3,0}x_{1,1}^2-G_{4,2}x_{1,0}, G_{4,2}x_{1,1}-G_{4,3}x_{1,0}, G_{4,2}x_{1,1}x_{1,0}^2-G_{4,3}G_{3,0},
G_{4,2}^2-G_{4,3}G_{3,0}x_{1,1}.\]

The advantage of this is that any monomial
\[x^ay^bz^c\mapsto (x_{1,1}x_{1,0}u)^a(x_{1,1})^b(x_{1,0}^3u^2)^c\]
can be reduced to normal form to determine where this lives.
For instance \[x^3y^2z\mapsto x_{1,1}^5x_{1,0}^6u^5\equiv G_{3,0}^2x_{1,1}^5
\equiv G_{4,2}^2x_{1,1}x_{1,0}^2\equiv G_{4,3}G_{4,2}G_{3,0}\]
shows this to be an element of $C(I^k,\ov{A})$ for any $k\leq 3$.


\newpage
\begin{thebibliography}{99}

\bibitem{dJ} 
 T. de Jong,
 {\emph An algorithm for computing the integral closure},
 J. Symbolic Computation, \textbf{26}, (1998), 273--277.

\bibitem{M2}
 D.R.Grayson and M.E.Stillman, 
 Macaulay2, a software system for research in algebraic geometry,
 Available at http://www.math.uiuc.edu/Macaulay2

\bibitem{ic}
 A.Taylor, D.Eisenbud, M.E.Stillman,
 IntegralClosure, a Macaulay2 package for computing integral closures of rings and ideals,
 Available at http://www.math.uiuc.edu/Macaulay2/Packages

%\bibitem{GP} 
% G.-M. Greuel and G. Pfister,
% {\tt normal.lib} A {\sc Singular} library for computing the normalization of affine rings, 2005.
%version="$id reesclos.lib,v 1.32 2001/01/16 hirsch Exp $";
%category="Commutative Algebra";

\bibitem{Hir}
T. Hirsch,
{\tt reesclos.lib} A {\sc Singular} library to compute the integral closure of an ideal, 2001.

%info="
%LIBRARY:     reesclos.lib   PROCEDURES TO COMPUTE THE INT. CLOSURE OF AN IDEAL
%AUTHOR:      Tobias Hirsch, email: hirsch@math.tu-cottbus.de

%OVERVIEW:
% A library to compute the integral closure of an ideal I in a polynomial ring
% R=k[x(1),...,x(n)] using the Rees Algebra R[It] of I. It computes the integral
% closure of R[It] (in the same manner as done in the library 'normal.lib'),
% which is a graded subalgebra of R[t]. The degree-k-component is the integral
% closure of the k-th power of I.

%PROCEDURES:
% ReesAlgebra(I);        computes the Rees Algebra of an ideal I
% normalI(I[,p[,r]]);    computes the integral closure of an ideal I using R[It]
%";

%LIB "normal.lib";       // for HomJJ
%LIB "standard.lib";     // for groebner

\bibitem{GPS}
 G.-M. Greuel, G. Pfister, and H. Sch\"onemann,
{\sc Singular} 3.1.5 A Computer Algebra System for Polynomial Computations,
Centre for Computer Algebra, University of Kaiserslautern, 2012.

\bibitem{qth} 
 D.A.Leonard,
 {\emph Finding the defining functions for one-point algebraic-geometric
  codes},
 IEEE Transactions on Information Theory, \textbf{47}, (2001), 2566--2573.  

\bibitem{Pell} 
 D.A.Leonard and R.Pellikaan, 
 {\emph Integral closures and weight functions over finite fields},
 Finite Fields and their Applications, \textbf{9}, (2003) 479--504.

\bibitem{mod}
 D.A.Leonard,
 {\emph A weighted module view of the integral closures of affine domains of type I},
 Advances in Mathematics of Communications, \textbf{3}, (2009), 1--11.

\bibitem{Qth}
D.A.Leonard,
 QthPower, a package to be available with Macaulay2 release 1.5 
 at http://www.math.uiuc.edu/Macaulay2/Packages

%\bibitem{Magma} 
% The {\sc MAGMA} 
% Computational Algebra System for Algebra, Number Theory and Geometry,
% The University of Sydney Computational Algebra Group.
% {\tt http://magma.maths.usyd.edu.au/magma}.

\bibitem{char0}
 D.A.Leonard,
 {\emph Extending the qth-power algorithm to integral closures of rings over the rationals},
 written in July, 2009, submitted in February, 2010.

\bibitem{Swan}
I. Swanson and C. Huneke,
{\emph Integral Closure of Ideals, Rings, and Modules},
London Mathematical Society Lecture Notes Series 336,
Cambridge University Press, 2006.

\bibitem{Vas06}
W. Vasconcelos,
{\emph Integral Closure: Rees Algebras, Multiplicities, Algorithms},
Springer Monographs in Mathematics,
Springer, 2005.

\bibitem{Vas98}
W. Vasconcelos,
{\emph Computational Methods in Commutative Algebra and Algebraic Geometry},
Algorithms and Computation in Mathematics 2,
Springer, 1998.

\end{thebibliography}
\end{document}